\documentclass[12pt,a4paper]{article}
\usepackage[utf8]{inputenc}
\usepackage[T2A]{fontenc}
\usepackage[english]{babel}
\usepackage{amsthm,amssymb, amsmath, amsfonts, amscd}
\usepackage{graphicx}
\usepackage{textcomp}
\usepackage{xcolor}
\usepackage{indentfirst}

\newtheorem{thm}{Theorem}
\newtheorem{prop}{Proposition}
\newtheorem{lm}{Lemma}
\newtheorem{rem}{Remark}
\newtheorem{cor}{Corollary}
\newtheorem{defin}{Definition}
\newtheorem*{defin-non}{Definition}

\def\beginproofdot{{\noindent\bf Proof. }}
\def\beginproof{{\noindent\bf Proof }}
\def\endproof{\hfill$\blacksquare$\medskip}

\begin{document}

\title{$L_2$-small ball asymptotics for a family of finite-dimensional perturbations of Gaussian functions}
\author{Petrova Yu.P.\footnote{Research is supported by «Native towns», a social investment program of PJSC «Gazprom Neft» and in part by the Moebius Contest Foundation for Young Scientists}}
\maketitle
\textit{Abstract:}
In this article we study the small ball probabilities in $L_2$-norm for a family of finite-dimensional perturbations of Gaussian functions. We define three types of perturbations: non-critical, partially critical and critical; and derive small ball asymptotics for the perturbated process in terms of the small ball asymptotics for the original process. The natural examples of such perturbations appear in statistics in the study of empirical processes with estimated parameters (the so-called Durbin's processes). We show that the Durbin's processes are critical perturbations of the Brownian bridge. Under some additional assumptions, general results can be simplified. As an example we find the exact $L_2$-small ball asymptotics for critical perturbations of the Green processes (the processes which covariance function is the Green function of the ordinary differential operator).

\section{Introduction}
The theory of small ball probabilities (also called small deviation probabilities) for various norms is extensively studied in recent decades (see the surveys~\cite{lif,lishao,2003fatalov}; for the extensive up-to-date bibliography see~\cite{Lifshits2018}). Given a random process $X(t)$, the asymptotic behavior of the probability $\mathbb{P}\left\{\|X\|<\varepsilon\right\}\text{ as } \varepsilon\to0$
is called an \textit{exact} asymptotics of small deviations. Let us note that in the literature on small balls, it is rare, and only known for a limited number of random processes, that the exact asymptotics is found. Therefore, often a \textit{ logarithmic asymptotics}  is studied, that is the asymptotics of  $\ln(\mathbb{P}\{\|X\|<\varepsilon\})$. But even on the logarithmic level the behavior of small balls cannot be uniformly described for the whole class of Gaussian measures.

In~\cite{sytaya} a solution of the small ball behavior problem was obtained in implicit terms. Starting from~\cite{ibragimov1979,zolotarev1979,hoffmann1979}, many authors attempted simplifying the asymptotic expression for the small ball probability under varios conditions.

Let $\mathcal{O}$ be a bounded domain in $\mathbb{R}^d$, $d\in\mathbb{N}$; $\bar{\mathcal{O}}$  is the closure of $\mathcal{O}$. Let $X_0(x)$, $x\in\bar{\mathcal{O}}$, be a zero mean-value  Gaussian random function with covariance function $G_0(x,y):=\mathbb{E}X_0(x)X_0(y)$. Denote by $\mathbb{G}_0$ the corresponding covariance operator in  $L_2(\mathcal{O})$:
\begin{align*}
(\mathbb{G}_0u)(x)=\int\limits_{\mathcal{O}}G_0(x,y)u(y)\,dy.
\end{align*}
Due to the Karhunen--Lo{\`e}ve expansion (see, e.g.~\cite[Section, 12]{lifshitz2016}) the small deviation problem can be rewritten as follows
\begin{align}\label{sb}
\mathbb{P}\left(\|X_0\|_{2}<\varepsilon\right):=\mathbb{P}\left(\int\limits_{\mathcal{O}}(X_0(x))^2\,dx<\varepsilon^2\right)=\mathbb{P}\left(\sum\limits_{k=1}^{\infty}\mu_k^0\xi_k^2<\varepsilon^2\right) \text{ as } \varepsilon\to0.
\end{align}

Here $\xi_k$, $k\in\mathbb{N}$, are the i.i.d. standard normal rv's, $\mu_k^0$ are the eigenvalues of the corresponding covariance operator $\mathbb{G}_0$. Moreover, due to the Karhunen--Lo{\`e}ve expansion we have $\sum_k\mu_k^0<\infty$. This condition means that $\|X_0\|_2$ is finite.

Therefore, knowing the eigenvalues $\mu_k^0$, one can obtain some information on the distribution of $\|X_0\|_2$. 
 In paper~\cite{DLL} Dunker, Lifshits and Linde found an exact asymptotics under some mild assumptions. In particular, the exact asymptotics was found in the cases $\mu_k^0=k^{-A}$  and $\mu_k^0=\exp(-Ak)$. The main difficulty is that we rarely know the explicit formulas for the eigenvalues. If one knows sufficiently precise asymptotics for $\mu_k^0$ then it is possible to obtain the small ball asymptotics up to a constant using well-known comparison principle of Wenbo Li:
\begin{prop}{\rm(\cite{Li1992,gao})} Let $\xi_k$ be a i.i.d. standard normal rv's; $\mu_k^0$ and $\mu_k$  be two positive non-increasing summable sequences such that  $\prod\mu_k/\mu_k^0<\infty$. Then
	\begin{align}\label{Li0}
	\mathbb{P}\Big\{\sum_{k=1}^{\infty} \mu_k^0 \xi_k^2<\varepsilon^2\Big\}\sim
	\mathbb{P}\Big\{\sum_{k=1}^{\infty} \mu_k \xi_k^2<\varepsilon^2\Big\}
	\cdot \left( \prod\limits_{k=1}^{\infty}\frac{\mu_k }{\mu_k^0} \right)^{1/2},
	\qquad \varepsilon\to0.
	\end{align}	
\end{prop}
In papers~\cite{naznik1,nazarov2009} there was selected the concept of the \textit{Green Gaussian process}. For such process the covariance function $G_0(x,y)$ is the Green function for an ordinary differential operator. This allows to study asymptotics for $\mu_k^0$ using the methods of spectral theory of ODEs, originated from the classical works of Birkhoff~\cite{birkhoff1,birkhoff2} and
Tamarkin~\cite{tamarkin1912,tamarkin1928} (further development of this theory can be found in~\cite{1982shkalikov,shkalikov1983}).
Many processes such as Wiener process, Brownian bridge, Ornstein-Uhlenbeck process and their integrated counterparts are Green processes.  Let us mention also an important series of works~\cite{2017KleptsynaGB,2017kleptsynafracBCO,2018KleptsynaFP}, where the exact asymptotics were obtained for non Green processes.


In this paper we consider a case of finite-dimensional perturbations of the Gaussian processes. It is known that the logarithmic asymptotics is not changing under finite-dimensional perturbation (in more general terms is proved in~\cite{2004GaoLi-loglevel} and~\cite{2009Naz-loglevel}). So we are interested in exact asymptotics. 

First the problem for one-dimensional perturbation was considered in \cite{deheuvels2007}. In paper~\cite{onedimpert} the notions of non-critical and critical perturbations were considered and exact small ball asymptotics were found for a wide class of one-dimensional perturbations.

In this article we generalize these results for the case of finite-dimensional perturbation  defined by the formula \eqref{Xpert}. 
 The motivating examples to consider such types of perturbations are the so-called Durbin's processes. These processes naturally appear in statistics as limiting ones when building goodness-of-fit tests of $\omega^2$-type for testing that a sample is belonging to the family of distributions with estimated parameters. 
 For example, if one considers test for normality (with estimated mean and/or variance) then  Kac-Kiefer-Wolfowitz processes (KKW processes) appear (see~\cite{KKW}). The exact $L_2$-small ball asymptotics for KKW processes is found in~\cite{SVFpert2015} and for some other important Durbin's processes in \cite{ZOO2017}.
 
The paper is organized as follows.  In Section \ref{finitedimpert} we describe the family of finite-dimensional perturbations of Gaussian random functions and define the notions of non-critical, partially critical and critical perturbations. In Section \ref{thm12} we prove genera1l theorems about $L_2$-small ball asymptotics in non-critical and critical cases (see Theorems \ref{thm1} and \ref{thm2}).   In Section \ref{Green} we apply these results in the case of the Green processes (Theorem 3). In Section \ref{ex-Durbin}  we consider the Durbin's processes and prove that they are the critical perturbations of the Brownian bridge (see Theorem \ref{Durbin-crit}). In Section \ref{sec:lemmas} (Appendix) we prove lemma about differentiability of the asymptotics of small ball probabilities. Besides the fact that we use this Lemma to prove Theorem 3, the Lemma is interesting by itself.

We use letter $C$ to denote various constants which exact values are not important. Denote by $E_m$ --- identical matrix of rank $m$.

\section{A set of finite-dimensional perturbations}
\label{finitedimpert}
%
Let $\vec{\varphi}(x)=(\varphi_1(x),\ldots,\varphi_m(x))^T$, where $\varphi_j(x)$ are locally summable functions in $\mathcal{O}$, $j=1\ldots m$.  Suppose that the vector-function
\begin{align*}
\vec{\psi}(x)=(\psi_1(x),\ldots,\psi_m(x))^T=\int\limits_{\mathcal{O}} G_0(x,y)\vec{\varphi}(y)\,dy
\end{align*}
is well-defined for a.e. $x\in\mathcal{O}$, $\psi_j\not\equiv0$, $j=1,\ldots, m$, and 
\begin{align}
Q_{ij}=\int\limits_{\mathcal{O}}\psi_i(x)\varphi_j(x)\,dx<\infty. \label{kappa}
\end{align}
This is equivalent to $\psi_j\in\mathrm{Im}(\mathbb{G}_0^{1/2})$.
\begin{rem}
	Without loss of generality we assume that the functions $\varphi_j(x)$, $j=1\ldots m$, are linearly independent.
\end{rem}
\begin{rem}\label{remQ}
	Formula~\eqref{kappa} defines the scalar product on the space dual to $\mathrm{Im}(\mathbb{G}_0^{1/2})$. Hence $Q$ is the Gram matrix, therefore, symmetric and nondegenerate.
\end{rem}   

We define a set of Gaussian functions
\begin{align}\label{Xpert}
X_{A}(x):=X_0(x)-\vec{\psi}(x)^T \cdot A\cdot\int\limits_{\mathcal{O}} X_0(y)\vec{\varphi}(y)\,dy.
\end{align}
 Here, $A$ is the matrix of perturbation parameters ($A_{ij}\in\mathbb{R}$, $i,j=1,\ldots,m$). Clearly, $\mathbb{E}X_{A}=0$. 
\begin{lm}
	 The covariance of $X_{A}(x)$ is
	\begin{align}\label{GXpert}
	G_{A}(x,y)
	=
	G_0(x,y)+\vec{\psi}(x)^T\cdot D\cdot\vec{\psi}(y),
	\end{align}
	where $D=-A-A^{T}+A QA^T.$ Note that $D$ is symmetric.
\end{lm}
\beginproofdot Formula \eqref{GXpert} can be checked by direct computation due to \eqref{kappa}.
\endproof
\begin{cor} For the processes \eqref{Xpert} the following equality holds
	\begin{align*}
	X_{A}(x)\stackrel{d}{=}X_{2Q^{-1}-A}(x).
	\end{align*}
\end{cor}

\begin{cor}\label{corcrit}
	Let $A=Q^{-1}$. Then the following assertions hold:
	\begin{enumerate}
		\item
			we have the identity a.s. {\normalfont{(}}$j=1,\ldots,m${\normalfont{)}}
			\begin{align*}
			\int\limits_{\mathcal{O}}X_{A}(x)\varphi_j(x)\,dx=0.
			\end{align*}	 
		\item
			the function $X_{A}(x)$ and the random variables $\int\limits_{\mathcal{O}} X_0(y)\varphi_j(y)\,dy$, $j=1,\ldots, m$, are independent.
		\item
			if $\varphi_j\in L_2(\mathcal{O})$, then the integral operator with the kernel function $G_{A}(x,y)$ has a zero eigenvalue of multiplicity $m$ corresponding to the eigenfunctions $\varphi_j$, $j=1,\ldots, m$.
	\end{enumerate}
\end{cor}
\beginproofdot
All statements follow from the relations:
\begin{align*}
\int\limits_{\mathcal{O}}X_{A}(x)\vec{\varphi}(x)^T\,dx
=\int\limits_{\mathcal{O}} X_0(x)\vec{\varphi}(x)^T\,dx\left(E_m-AQ\right);\\
\mathbb{E}X_{A}(x)\int\limits_{\mathcal{O}} X_0(y)\vec{\varphi}(y)^T\,dy
=
\vec{\psi}(x)^T(E_m-AQ),
\end{align*}
which can be easily checked.
\endproof
\begin{defin}\label{def:nocrit}
	A function $X_{A}$ is called a non-critical perturbation of the function $X_0$, if the following equivalent conditions hold:
	\begin{enumerate}
		\item  $\det(E_m-A^T Q)\not=0$;
		\item 
		$\int\limits_{\mathcal{O}}X_{A}(x)\varphi_j(x)\,dx$, $j=1,\ldots,m$, are linearly independent a.s.
	\end{enumerate} 
\end{defin}
\begin{defin}\label{def:semicrit}
		A function $X_{A}$ is called a partially critical perturbation of rank $s$ of the function $X_0$, $0< s< m$, if the following equivalent conditions hold:
	\begin{enumerate}
		\item  $\mathrm{rank}(E_m-A^T Q)=m-s$;
		\item 
		$\int\limits_{\mathcal{O}}X_{A}(x)\varphi_j(x)\,dx$, $j=1,\ldots,m$, form a vector space of dimension $m-s$ a.s.
	\end{enumerate} 
\end{defin}
\begin{defin}\label{def-crit}
		A function  $X_{A}$ is called a critical perturbation of the function $X_0$, if the following equivalent conditions hold:
		\begin{enumerate}
			\item  $A=Q^{-1}$;
			\item 
			$\int\limits_{\mathcal{O}}X_{A}(x)\varphi_j(x)\,dx=0$, $j=1,\ldots,m$ a.s.
		\end{enumerate} 
\end{defin}
\begin{rem}
	For one-dimensional case the notions of critical and non-critical perturbations were introduced in \cite{onedimpert}.
\end{rem}

\section{Main theorems (non-critical  and critical case)}
\label{thm12}
Let $\mu_k^0$ и $u_k^0(x)$ be the eigenvalues and the corresponding (normalized in $L_2(\mathcal{O})$) eigenfunctions of the covariance operator $\mathbb{G}_{0}$, that is 
\begin{align*}
\mu_k^0 u_k^0(x)=\int\limits_{\mathcal{O}}G_{0}(x,y)u_k^0(y)\,dy.
\end{align*}
As $\|X_0\|_2<\infty$ a.s., the operator $\mathbb{G}_0$ belongs to the trace class, that is $\sum_k\mu_k^0<\infty$. 
 Let $\mu_k$ and $u_k(x)$ be the eigenvalues and the corresponding eigenfunctions of the integral operator with kernel $G_{A}(x,y)$. Note that due to the minimax principle (see, e.g., \cite[Section 10.2]{birmansolomENG}), the sequences $\mu_k^0$ and $\mu_k$ interlace. In particular this implies the convergence of the series $\sum_k\mu_k$ and $\|X_A\|_2<\infty$ a.s.
 By definition, put $\lambda_k^0:=(\mu_k^0)^{-1}$, $\lambda_k:=\mu_k^{-1}$. 

\begin{thm}\label{thm1}	Let $X_A$ be a non-critical perturbation of function $X_0$. As $\varepsilon\to0$ we have
	\begin{align*}
		\mathbb{P}\left\{\|X_{A}\|_2<\varepsilon\right\}\sim \frac{\mathbb{P}\left\{\|X_0\|_2<\varepsilon\right\}}{\det\left(E_m-QA\right)}.
	\end{align*}
\end{thm}
\beginproofdot
By the comparison theorem (see Prop.~\ref{Li0}) as $\varepsilon\to0$ we obtain
\begin{align*}
\mathbb{P}\left\{\|X_{A}\|_2<\varepsilon\right\}\sim \mathbb{P}\left\{\|X_0\|_2<\varepsilon\right\}\cdot \left(\prod\limits_{k=1}^{\infty}\frac{\mu_{k}^0}{\mu_k}\right)^{1/2}.
\end{align*}

Consider the Fredholm determinants for the kernels $G_0$ and $G_{A}$, respectively:
\begin{align*}
\mathcal{F}^0(z)
:=\prod\limits_{k=1}^{\infty}\left(1-\frac{z}{\lambda_k^0}\right);
\qquad 
\mathcal{F}(z)
:=\prod\limits_{k=1}^{\infty}\left(1-\frac{z}{\lambda_k}\right).
\end{align*}
Since the series $\sum\limits_k (\lambda_k^0)^{-1}$ and $\sum\limits_k \lambda_k^{-1}$ converge, these canonical Hadamard products converge for all $z\in\mathbb{C}$. Jensen's theorem (see \cite[\S3.6]{titch}) provides
\begin{align}\label{Iensen1}
\prod\limits_{k=1}^{\infty}\frac{\lambda_k^0}{\lambda_{k}}
=
\lim\limits_{|z|\to\infty}\exp\left(\frac{1}{2\pi}\int\limits_0^{2\pi}\ln\left|\frac{\mathcal{F}(z)}{\mathcal{F}^0(z)}\right|\,d\arg(z)\right).
\end{align}
Using the transformation formula for Fredholm determinant under finite-dimensional perturbation (see~\cite{1908bateman},~\cite[Theorem  2.2]{1972sukhatme}, or more general formula~\cite[Chapter 2, \S 4.6]{1950kantorovichkrylov}) we obtain
\begin{align}\label{Sukhatme}
\frac{\mathcal{F}(z)}{\mathcal{F}^0(z)}=\det (L(z)),
\end{align}
where matrix $L(z)$ is defined from the formula
\begin{align}\label{L}
	L(z)=E_m+\sum\limits_{n=1}^{\infty}\frac{\lambda_n^0\vec{a}_n\vec{a}_n^{\,T}}{1-\frac{\lambda_n^0}{z}}\cdot D.
\end{align}
Here, $\vec{a}_n=\int\limits_{\mathcal{O}}\vec{\psi}(x)u_n^0(x)\,dx$ are the Fourier coefficients of vector-function $\vec{\psi}(x)$ with respect to the orthogonal basis $\{u_k^0\}$.

To justify the passage to the limit in~\eqref{Iensen1} we argue similarly to~\cite[Lemma 5.1]{onedimpert}. As $|z|\to\infty$ the entries of matrix $L(z)$ converge uniformly with respect to $\arg(z)\in[\varepsilon,2\pi-\varepsilon]$, $\varepsilon>0$ to the entries of matrix $\left(E_m+\sum\limits_{n=1}^{\infty}\lambda_n^0\vec{a}_n\vec{a}_n^{\;T}\cdot D\right)$. In the neighborhood of positive real line the integrand has a summable majorant, hence by Lebesgue dominated convergence theorem the limit~\eqref{Iensen1} is equal
\begin{align}
\label{S1:det1}
\prod\limits_{k=1}^{\infty}\frac{\lambda_k^0}{\lambda_{k}}
=\det\left(E_m+\sum\limits_{n=1}^{\infty}\lambda_n^0\vec{a}_n\vec{a}_n^{\;T}\cdot D\right).
\end{align}
Note that 
\begin{align*}
\vec{\psi}(x)
=\sum\limits_{n=1}^{\infty}\vec{a}_n u_n^0(x);
&&
\vec{\varphi}(x)
=\sum\limits_{n=1}^{\infty}\lambda_n^0\vec{a}_n u_n^0(x).
\end{align*}
Then by the orthonormality of $u_n^0$
\begin{align}\label{S1:Q}
Q
&=\int\limits_{\mathcal{O}}\vec{\varphi}(x)\vec{\psi}^{\;T}(x)\,dx
=\sum\limits_{n=1}^{\infty}\lambda_n^0\vec{a}_n\vec{a}_n^{\;T}.
\end{align}
Hence the formula~\eqref{S1:det1} has the following form:
\begin{align*}
\prod\limits_{k=1}^{\infty}\frac{\lambda_k^0}{\lambda_{k}}
&=\det\left(E_m+Q\cdot D\right)
=
\det\left(E_m+Q\cdot[-A-A^{T}+A^T\cdot Q\cdot A)]\right)
\\
&=
\det\left(E_m-Q\cdot A^T\right)\cdot \det\left(E_m-Q\cdot A\right)
\\
&=
\det\left(E_m-A\cdot Q\right)\cdot \det\left(E_m-Q\cdot A\right)=\left(\det\left(E_m-Q\cdot A\right)\right)^2.
\end{align*}
The last equality holds due to similarity of matrices $E_m-QA$ and $E_m-AQ$. It's clear that in non-critical case $\det\left(E_m-QA\right)\not=0$. This completes the proof.
\endproof


To work with critical perturbations we need the following condition.\\
\textbf{Condition A:} $\forall j=1,\ldots, m:\quad	\varphi_j\in L_2(\mathcal{O}),\;\text{that is equivalent to } \; \psi_j\in\mathrm{Im}(\mathbb{G}_0).
$

\begin{thm}
	\label{thm2}
	 Let $X_{A}$ be a critical perturbation of function $X_0$ (see Definition 3) and the condition A be fulfilled.
	 Then as $r\to0$ 
	\begin{align}\label{thm2Pasymp}
	\nonumber
	&\mathbb{P}\left\{\|X_{A}\|_2<\sqrt{r}\right\}
	\sim\sqrt{ \frac{\det(Q)}{\det\left(\int\limits_{\mathcal{O}}\vec{\varphi}(s)\vec{\varphi}^{\;T}(s)\,ds\right)}}\left(\sqrt{\frac{2}{\pi}}\right)^{m}
	\\
	&\times \int\limits_0^r\int\limits_0^{r_1}\ldots \int\limits_0^{r_{m-1}}\frac{d^m}{dr_m^m}\mathbb{P}\bigl\{\|X_0\|_2<\sqrt{r_m}\bigr\} \frac{dr_{m}\ldots dr_1}{\sqrt{(r-r_1)(r_1-r_2)\ldots(r_{m-1}-r_{m})}}.
	\end{align}
\end{thm}
\beginproofdot 
We introduce three distribution functions:
\begin{align*}
F_0(r):
&=\mathbb{P}\Bigl\{\sum\limits_{k=1}^{\infty}\mu_k^0\xi_k^2<r\Bigr\}=\mathbb{P}\Bigl\{\|X_0\|_2<\sqrt{r}\Bigr\};
\\
F(r):
&=\mathbb{P}\Bigl\{\sum\limits_{k=1}^{\infty}\mu_k\xi_k^2<r\Bigr\}=\mathbb{P}\Bigl\{\|X_{A}\|_2<\sqrt{r}\Bigr\};
\\
F_m(r):
&=\mathbb{P}\Bigl\{\sum\limits_{k=m+1}^{\infty}\mu_k^0\xi_k^2<r\Bigr\},\quad m\in\mathbb{N}.
\end{align*}
Let us show that as $r\to0$ 
\begin{align}\label{FFm}
F(r)\sim F_m(r)\cdot \left(\prod\limits_{k=1}^{\infty}\frac{\mu_{k+m}^0}{\mu_{k}}\right)^{1/2}.
\end{align}
The Jensen's theorem provides
\begin{align}\label{Iensen}
\prod\limits_{k=1}^{\infty}\frac{\mu_{k}}{\mu_{k+m}^0}
=\prod\limits_{k=1}^{\infty}\frac{\lambda_{k+m}^0}{\lambda_{k}}=
\lim\limits_{|z|\to\infty}\exp\left(\frac{1}{2\pi}
\int\limits_0^{2\pi}\ln\left|\frac{\mathcal{F}(z)}{\mathcal{F}^0(z)}\cdot\prod\limits_{l=1}^{m}\left(1-\frac{z}{\lambda_l^0}\right)\right|\,d\arg(z)\right).
\end{align}
Note that in critical case $E_m=-QD$. Thus using~\eqref{S1:Q} the matrix~\eqref{L} can be rewritten as follows:
\begin{align}\label{L1}
L&=-\sum\limits_{n=1}^{\infty}\lambda_n^0\vec{a}_n\vec{a}^{\;T}_n\cdot D+\sum\limits_{n=1}^{\infty}\frac{\lambda_n^0\vec{a}_n\vec{a}^{\;T}_n}{1-\frac{\lambda_n^0}{z}}\cdot D=\frac1z\sum\limits_{n=1}^{\infty}\frac{(\lambda_n^0)^2\vec{a}_n\vec{a}_n^{\;T}}{1-\frac{\lambda_n^0}{z}}\cdot D.
\end{align}
Hence, due to~\eqref{Sukhatme} and~\eqref{L1}, the expression under the log sign in formula~\eqref{Iensen} can be rewritten as follows: \begin{align*}
\left|\frac{\mathcal{F}(z)}{\mathcal{F}^0(z)}\cdot\prod\limits_{l=1}^{m}\left(1-\frac{z}{\lambda_l^0}\right)\right|
&=
\left|\det(L_{ij})_{i,j=1}^m\cdot\prod\limits_{l=1}^{m}\left(1-\frac{z}{\lambda_l^0}\right)\right|
\\
&=
\left|\det\sum\limits_{n=1}^{\infty}\frac{(\lambda_n^0)^2\vec{a}_n\vec{a}_n^{\;T}}{1-\frac{\lambda_n^0}{z}}\cdot D\cdot\prod\limits_{l=1}^{m}\left(\frac{1}{z}-\frac{1}{\lambda_l^0}\right)\right|.
\end{align*}
By~\cite[Lemma 5.1]{onedimpert} the limit in right hand side of formula~\eqref{Iensen} equals: 
\begin{align}\label{eigenprod}
\prod\limits_{k=1}^{\infty}\frac{\mu_{k}}{\mu_{k+m}^0}
=
\left|\det\left(\sum\limits_{n=1}^{\infty}(\lambda_n^0)^2\vec{a}_n\vec{a}_n^{\;T}\cdot D\right)\right|\cdot\prod\limits_{l=1}^{m}\frac{1}{\lambda_l^0}.
\end{align}

Consider the well-defined matrix 
\begin{align*}
	\int\limits_{\mathcal{O}}\vec{\varphi}(x)\vec{\varphi}^{\;T}(x)\,dx
	&=
	\int\limits_{\mathcal{O}}
	\left(\sum\limits_{k=1}^{\infty}\lambda_k^0\vec{a}_k  u_k(x)\right)
	\left(\sum\limits_{n=1}^{\infty}\lambda_n^0\vec{a}_n u_n(x)\right)^T
\,dx
	=\sum\limits_{n=1}^{\infty}(\lambda_n^0)^2\vec{a}_n\vec{a}_n^{\;T}.
\end{align*}
Then
\begin{align*}
\prod\limits_{k=1}^{\infty}\frac{\mu_{k}}{\mu_{k+m}^0}
=\det\left(\int\limits_{\mathcal{O}}\vec{\varphi}(x)\vec{\varphi}^{\;T}(x)\,dx\right)\cdot\frac{1}{\det(Q)}\cdot\prod\limits_{l=1}^{m}\frac{1}{\lambda_l^0}.
\end{align*}

Combining formulas~\eqref{FFm} and~\eqref{eigenprod}, we finally get as $r\to0$
\begin{align}\label{FtildeFk}
F(r)\sim F_m(r)\cdot\sqrt{ \frac{\det(Q)\cdot\lambda_1^0\cdot\ldots\cdot\lambda_m^0}{\det\left(\int\limits_{\mathcal{O}}\vec{\varphi}(x)\vec{\varphi}^{\;T}(x)\,dx\right)}}.
\end{align}

Further, obviously,
$F_0(r_m)=(F_m*f_1*\ldots*f_m)(r_m)$, where
\begin{align*}
f_j(x)
=
\frac{d}{dx}\mathbb{P}\{\mu_j^0\xi_j^2\leqslant x\}
=
\begin{cases}
\frac{\exp\left(-\frac{x}{2\mu_j^0}\right)}{\sqrt{2\pi\mu_j^0 x}},&x> 0,
\\
0,&x\leqslant0,
\end{cases}
\qquad
 j=1\ldots m.
\end{align*}

Note that the following relation holds
\begin{align*}
(F_{m}*f_{m})(r)=F_{m-1}(r).
\end{align*}
By the Laplace transform we obtain the solution of this convolution equation:
\begin{align*}
F_{m}(z)
&=\sqrt{\frac{2\mu^0_{m}}{\pi}}\int\limits_0^z\left(F'_{m-1}(r_1)+\frac1{2\mu_{m}^0}F_{m-1}(r_1)\right)\exp\left(-\frac{z-r_1}{2\mu_{m}^0}\right)\frac{dr_1}{\sqrt{z-r_1}}.
\end{align*}
By Lemma~\ref{lemmaF} from Appendix we obtain $F_{m-1}(r_1)=o(F_{m-1}'(r_1))$, $r_1\to+0$. Hence, as $z\to+0$
\begin{align*}
F_{m}(z)
&\sim\sqrt{\frac{2\mu_{m}^0}{\pi}}\int\limits_0^zF'_{m-1}(r_1)\frac{dr_1}{\sqrt{z-r_1}}.
\end{align*}
 Analogously, due to relation between $F_{m-1}$ and $F_{m-2}$, we obtain
\begin{align*}
F_{m-1}(r_1)
&=\sqrt{\frac{2\mu_{m-1}^0}{\pi}}\int\limits_0^{r_1}\left(F'_{m-2}(r_2)+\frac1{2\mu_{m-1}^0}F_{m-2}(r_2)\right)\exp\left(-\frac{r_1-r_2}{2\mu_{m-1}^0}\right)\frac{dr_2}{\sqrt{r_1-r_2}};
\\
F_{m-1}'(r_1)
&=\sqrt{\frac{2\mu_{m-1}^0}{\pi}}\int\limits_0^{r_1}\left(F''_{m-2}(r_2)+\frac1{2\mu_{m-1}^0}F'_{m-2}(r_2)\right)\exp\left(-\frac{r_1-r_2}{2\mu_{m-1}^0}\right)\frac{dr_2}{\sqrt{r_1-r_2}}.
\end{align*}
By Lemma~\ref{lemmaF} from Appendix we have $F_{m-1}'(r_2)=o(F_{m-1}''(r_2))$, $r_2\to+0$. Hence, as $r_1\to+0$ we obtain
\begin{align*}
F_{m-1}'(r_1)
&\sim\sqrt{\frac{2\mu_{m-1}^0}{\pi}}\int\limits_0^{r_1}F''_{m-2}(r_2)\frac{dr_2}{\sqrt{r_1-r_2}}.
\end{align*}
Thus
\begin{align*}
F_{m}(z)
&\sim\left(\sqrt{\frac{2}{\pi}}\right)^2\sqrt{\mu_{m}^0\mu_{m-1}^0}\int\limits_0^z \int\limits_0^{r_1}F''_{m-2}(r_2)\frac{dr_2\,dr_1}{\sqrt{(z-r_1)(r_1-r_2)}}.
\end{align*}
Using this algorithm $m-2$ times,  we get  as $z\to +0$
\begin{align}\label{Fk+1asymp}
F_{m}(z)
&\sim\left(\sqrt{\frac{2}{\pi}}\right)^{m}\prod\limits_{l=1}^{m}{\sqrt{\mu_{l}^0}}\int\limits_0^z\int\limits_0^{r_1}\ldots \int\limits_0^{r_{m-1}} \frac{F_0^{(m)}(r_{m})\,dr_{m}\ldots dr_1}{\sqrt{(z-r_1)(r_1-r_2)\ldots(r_{m-1}-r_{m})}}.
\end{align}
Combining the relations~\eqref{FtildeFk} and~\eqref{Fk+1asymp} completes the proof.
\endproof
\begin{rem}
	The results analogous to Theorems \ref{thm1} and \ref{thm2} for one-dimensional perturbations were obtained in~\cite{onedimpert}.
\end{rem}
\begin{rem}
	In case of partially critical perturbation  {\normalfont{(}}see Definition~\ref{def:semicrit}{\normalfont{)}}, if the condition A is fulfilled, then the asymptotics of small ball probabilities can be found using Theorems~\ref{thm1} and~\ref{thm2}.
\end{rem}
\section{Small ball probabilities for Green processes}
\label{Green}
Now we suppose that $\mathcal{O}=(0,1)$, and the covariance function $G_0(t,s)$, $t,s\in[0,1]$, is the Green function of a self-adjoint operator $L_0$ in the space $L_2(0,1)$, generated by a differential expression of order $2l$:
\begin{align}\label{L0}
L_0u:=(-1)^lu^{(2l)}+(p_{l-1}u^{(l-1)})^{(l-1)}+\ldots+p_0u,
\end{align}
and $2l$ boundary conditions. We recall that by definition, $G_0$ for any $s\in(0,1)$ satisfies the equation $L_0G_0=\delta(t-s)$, in the sense of distributions, and satisfies boundary conditions.

By $\mathcal{D}(L_0)$ we denote the image of an integral operator with the kernel function $G_0(s,t)$. Then the inverse operator is just $L_0$ with the domain $\mathcal{D}(L_0)$. In particular, if $\varphi\in L_2(0,1)$, then $\psi\in\mathcal{D}(L_0)$, and $L_0\psi=\varphi$. Assume for simplicity that $p_j\in C^j[0,1]$. Then $\mathcal{D}(L_0)$ coincides with the set of functions which belong to $W_2^l(0,1)$ and satisfy boundary conditions.

\begin{thm}
	
	Let $A=Q^{-1}$. If 
	$\varphi_i(x)\in L_2(\mathcal{O})$, $i=1,\ldots,m$, then as $\varepsilon\to0$
	\begin{align}\label{Green2}
	\mathbb{P}\{\|X_A\|_2\leqslant\varepsilon\}\sim \sqrt{\frac{\det(Q)}{\det\left(\int\limits_{\mathcal{O}}\vec{\varphi}(s)\vec{\varphi}^{\;T}(s)\,ds\right)}}\cdot (2l\sin(\pi/(2l))\varepsilon^2)^{-\frac{lm}{2l-1}}\cdot \mathbb{P}\{\|X_0\|_2\leqslant\varepsilon\}.
	\end{align}
\end{thm}
\beginproofdot
It is shown in \cite{naznik1} that the process $X_0$ satisfies the relation
\begin{align}\label{Fasymp}
F(r)=\mathbb{P}\{\|X_0\|_2\leqslant\sqrt{r}\}\sim C\cdot r^{\beta}\exp(-\mathcal{D}r^{-d}),\qquad r\to0,
\end{align}
where $d=\frac{1}{2l-1}$, $\mathcal{D}=\frac{1}{2d}(2l\sin(\pi/(2l)))^{-d-1}$.
Using Lemma \ref{diff} (see Appendix),  we obtain
 \begin{align}\label{Fmasymp}
 F^{(m)}(r)\sim C\cdot(\mathcal{D}d)^m
\cdot r^{\beta-m(d+1)} \exp(-\mathcal{D}r^{-d}), \qquad r\to0.
\end{align}
This means that the asymptotics \eqref{Fasymp} is $m$-times differentiable with respect to $r$. Substituting \eqref{Fmasymp} into \eqref{thm2Pasymp}, we obtain
\begin{align*}
	\mathbb{P}\{\|X_A\|_2\leqslant \sqrt{r}\}
	&\sim\sqrt{ \frac{\det(Q)}{\det\left(\int\limits_{\mathcal{O}}\vec{\varphi}(s)\vec{\varphi}^{\;T}(s)\,ds\right)}}\left(\sqrt{\frac{2}{\pi}}\right)^{m}\cdot C (\mathcal{D}d)^m
	\\
	&\times \int\limits_0^r\int\limits_0^{r_1}\ldots \int\limits_0^{r_{m-1}} \frac{r_m^{\beta-m(d+1)} \exp(-\mathcal{D}r_m^{-d}) \;dr_{m}\ldots dr_1}{\sqrt{(r-r_1)(r_1-r_2)\ldots(r_{m-1}-r_{m})}}.
\end{align*}
Consider the integral
\begin{align*}
I=\int\limits_0^{r_{m-1}} \frac{r_m^{\beta-m(d+1)} \exp(-\mathcal{D}r_m^{-d}) \;dr_{m}}{\sqrt{r_{m-1}-r_{m}}}.
\end{align*}
By \cite[Theorem 3]{onedimpert}, we obtain
\begin{align*}
I\sim \sqrt{\frac{\pi}{\mathcal{D}d}}\cdot r_{m-1}^{\beta-m(d+1)+(d+1)/2}\exp(-\mathcal{D}r_{m-1}^{-d}), \qquad r_{m-1}\to0.
\end{align*}
Repeating this procedure $(m-1)$ times completes the proof.
\endproof
\section{Example: Durbin's processes}
\label{ex-Durbin}

Durbin's processes appear as limit ones when building goodness-of-fit tests of $\omega^2$-type for testing if a sample belongs to a family of distributions when parameters of the distribution are estimated from the sample (see \cite{durbin1973}). For reader's convenience we describe them here.

Let $x_1, \ldots, x_n\in\mathbb{R}$ be a sample from a distribution with general distribution function $F(x,\theta)$, let $f(x,\theta)$ be the distribution density, where $\theta=(\theta_1,\ldots,\theta_s)$, $s\in\mathbb{N}$, is a vector of parameters.  Consider an empirical distribution function for fixed parameter values  $\theta^{\,0}=(\theta_1^0,\ldots,\theta_s^0)$: 
$$F_{n}^{0}(t)=\frac{\#\{x_i \,\colon F(x_i,\theta^{\,0})\leqslant t,\,i=1,\ldots,n\}}{n}, \quad t\in[0,1].$$

It is well known (see~\cite[Chapter 3]{Bil1977}), that the process  $n^{1/2}\bigl[F_{n}^0(t)-t\bigr]$ weakly converges to the Brownian bridge $B(t)$ in Skorokhod space $D[0,1]$.

Assume that some parameters of the distribution are unknown (without loss of generality we may assume that these are the first $m$ parameters). The unknown parameters are estimated using the sample (e.g., by applying the maximum likelihood method), and the new parameter vector is denoted by  $\hat{\theta}:=(\hat{\theta}_1,\ldots,\hat{\theta}_m,\theta_{m+1}^0,\ldots,\theta_s^0)$.
Then the empirical distribution function becomes
$$
\hat{F}_{n}(t)=\frac{\#\{ x_i\colon F(x_i,\hat{\theta})\leqslant t,\,i=1,\ldots,n\}}{n}, \quad t\in[0,1].
$$
It was shown in~\cite{durbin1973} that the process $n^{1/2}\bigl[\hat{F}_{n}(t)-t\bigr]$ converges weakly in $D[0,1]$ to a finite-dimensional perturbation of the Brownian bridge, namely, to a Gaussian process with zero mean and the covariance function
\begin{align}\label{Durbincov}
	G(s,t)=G_B(s,t)-\vec{\psi}^{\;T}(s)\,S^{-1}\,\vec{\psi}(t),\qquad s,t\in[0,1].
\end{align}
Here,	$G_B(s,t)=\min(s,t)-st$ is the covariance function of the Brownian bridge $B(t)$, $S$ is the Fisher information matrix with entries $S_{ij}$, $i,j=1,\ldots,m$:
\begin{align*}
S_{ij}=-\mathbb{E}\,\Bigl(\frac{\partial^2 }{\partial \theta_i\partial\theta_j} \ln(f(x,\theta))\Bigr)\biggr\rvert_{\theta=\theta_0}
=
\mathbb{E}\,\Bigl(\frac{\partial }{\partial \theta_i}\ln(f(x,\theta))\frac{\partial}{\partial\theta_j} \ln(f(x,\theta))\Bigr)\biggr\rvert_{\theta=\theta_0}.
\end{align*}
And the vector function  $\vec{\psi}=\bigl(\psi_1(t),\ldots,\psi_m(t)\bigr)$ is defined as 
\begin{align*}
\psi_j(t)=\frac{\partial F(x,\theta)}{\partial \theta_j}\Bigr\rvert_{\theta=\theta^0},\quad j=1,\ldots,m,
\end{align*}
where  $x$ and $t$ are related by $t=F(x,\theta)$.
Formula \eqref{Durbincov} shows that Durbin's processes are the $m$-dimensional perturbations of the Brownian bridge (of type~\eqref{Xpert}).
\begin{rem} The following equality holds:
\begin{align}\label{q'}
\psi_j'(t)
&=\frac{\partial }{\partial \theta_j}\ln(f(F^{-1}(t,\theta)))\biggr\rvert_{\theta=\theta_0}.
\end{align}
\end{rem}
Formula \eqref{q'} can be checked by direct computation. 
Thus from~\eqref{q'} we obtain
\begin{align*}
S_{ij}=\int\limits_0^1 \psi_i'(t)\psi_j'(t)\,dt.
\end{align*}

\begin{thm}\label{Durbin-crit}
	The Durbin's processes with $m$ estimated parameters are critical. 
\end{thm}
\beginproofdot
Note that if $X$ is the Brownian bridge, then we have $\varphi_i(s)=-\psi_i''(s)$, and
\begin{align*}
Q_{ij}=\int\limits_0^1\psi_j(s)\varphi_i(s)\,ds=\int\limits_0^1\psi_j(s)(-\psi_i''(s))\,ds=\int\limits_0^1\psi_j'(s)\psi_i'(s)\,ds=S_{ij}.
\end{align*}
Hence from~\eqref{def-crit} and~\eqref{Durbincov} the statement of the theorem follows.
\endproof

\section{Appendix}
Our aim is to prove the following lemma on differentiability of asymptotics of small ball probabilities.
\label{sec:lemmas}
\begin{lm}\label{diff}
Let $$F(x)=\mathbb{P}\Bigl\{\sum\limits_{k=1}^{\infty}\mu_k\xi_k^2<x\Bigr\},$$
where $\mu_k>0$, $k\in\mathbb{N}$, $\sum\limits_{k=1}^{\infty}\mu_k<\infty$, and $\xi_k$ are i.i.d. standard normal random variables.\\	
If $F(x)$ have the following asymptotics at zero
\begin{align*}
	F(x)\sim \mathcal{A}\,x^{\alpha}L(x)\exp(-\mathcal{D}x^{-\beta}),\quad \text{ as }x\to+0,\qquad \alpha\in\mathbb{R}, C>0,\beta>0, \mathcal{D}>0,
\end{align*}
where $L(x)>0$ is a slowly varying function at zero. \\Then for any $m\in\mathbb{N}$ as $x\to+0$
\begin{align*}
F^{(m)}(x)\sim \mathcal{A}(\mathcal{D}\beta)^{m}\,x^{\alpha-m(\beta+1)}L(x)\exp(-\mathcal{D}x^{-\beta}).
\end{align*}
\end{lm}
The proof of Lemma \ref{diff} is based on some good properties of $F(x)$ (Lemma \ref{lemmaF}) and the lemma of Tauberian type (Lemma \ref{Tauberian}).
\begin{lm}\label{lemmaF}
Let $F(x)$ satisfy the conditions of Lemma \ref{diff}. \\
Then $F^{(n)}(0)=0,\; n\in\mathbb{N}\cup\{0\}.$ And  $F^{(n)}(x)=o(F^{(n+1)}(x)),\quad x\to+0.$

\end{lm}
\beginproofdot
\textbf{\textit{Step 1:}} Assume that the following relations hold in a neighborhood of $x=0$, $x>0$, 
\begin{align}
	\label{step0}
	F^{(n+2)}(x)>0;\qquad F^{(n+1)}(x) \text{ is bounded\quad and\quad } F^{(n)}(0)=0.
\end{align}   
We claim that  $F^{(n)}(x)=o(F^{(n+1)}(x))$, $x\to+0$.
Indeed, integrating by parts we get
\begin{align*}
	\int\limits_0^x y F^{(n+2)}(y)\,dy=y F^{(n+1)}(y)\biggr\rvert_0^x-\int\limits_0^x F^{(n+1)}(y)\,dy=x F^{(n+1)}(x)-F^{(n)}(x).
\end{align*}
If we take $x$ so small that $F^{(n+2)}(y)>0$ for all $y\in(0,x)$, then the integral is positive. That means $x F^{(n+1)}(x)-F^{(n)}(x)>0$, and $x F^{(n+1)}(x)> F^{(n)}(x)$. This results in  $F^{(n)}(x)=o(F^{(n+1)}(x))$, $x\to+0$.

\textbf{\textit{Step 2:}} Let us prove Lemma~\ref{lemmaF} for the distribution function of the  \emph{finite} sum $\eta_m:=\sum\limits_{j=1}^{m}\mu_j\xi_j^2$. Namely, consider
\begin{align*}
	F_{\eta_m}(x):=\mathbb{P}\left\{\eta_m<x\right\}=\mathbb{P}\Bigl\{\sum\limits_{j=1}^{m}\mu_j\xi_j^2<x\Bigr\}.
\end{align*}
We claim that for all natural  $n<m/2-1$ the following relations hold
\begin{align}
	\label{step1}
	F_{\eta_m}^{(n+1)}(0)=0\qquad\text{and}\qquad
	F^{(n)}_{\eta_m}(x)=o(F^{(n+1)}_{\eta_m}(x)),\qquad x\to+0.
\end{align}
Changing variables in $F_{\eta_m}(x)$ to spherical coordinates, we obtain:
\begin{align*}
	F_{\eta_m}(x)=\int\limits_0^{\sqrt{x}}\int\limits_{S^{m-1}}d\varphi_1\ldots d\varphi_{m-1}\cdot  e^{-r^2\cdot P(\sin(\varphi_1),\ldots,\sin(\varphi_{m-1}))}\cdot |J|\,\frac{dr}{\sqrt{\mu_1\cdot\ldots\cdot\mu_m}},
\end{align*}
where $P(y_1,\ldots,y_{m-1})$ is a polynomial, $S^{m-1}$ is a $(m-1)$-dimensional sphere in $\mathbb{R}^m$, $J$  is the Jacobian:
\begin{align*}
	J=r^{m-1}\cdot \sin(\varphi_2)\cdot\sin^2(\varphi_3)\cdot\ldots\cdot\sin^{m-2}(\varphi_{m-1}).
\end{align*}
\color{red}
\color{black}\\
Hence, the density  $f_{\eta_m}(x)=F_{\eta_m}'(x)$  is equal to
\begin{align}
	\label{fetam}
	f_{\eta_m}(x)=&x^{m/2-1}\cdot\frac12\cdot\int\limits_{S^{m-1}}d\varphi_1\ldots d\varphi_{m-1}\cdot  e^{-x\cdot P(\sin(\varphi_1),\ldots,\sin(\varphi_{m-1}))}
	\\
	&\times|\sin(\varphi_2)\cdot\ldots\cdot\sin^{m-2}(\varphi_{m-1})|
	\,\frac{1}{\sqrt{\mu_1\cdot\ldots\cdot\mu_m}}.
	\nonumber 
\end{align}
It follows from~\eqref{fetam} that while $n<m/2-1$ the derivative $f_{\eta_m}^{(n)}(x)$ is defined in the neighborhood of $x=0$, and $f_{\eta_m}^{(n)}(0)=0$. Moreover,  $f_{\eta_m}^{(n)}(x)\sim C\cdot x^{m/2-1-n}$ as $x\to+0$ for some constant $C>0$, hence,~\eqref{step1} holds. Note that while $n<m/2-1$ we also have $F_{\eta_m}^{(n+2)}(x)>0$ in some neighborhood of $x=0$.


\textbf{\textit{Step 3:}} Let $\eta$ be a random variable (independent of  $\xi_1, \ldots, \xi_m$)  on the positive  half-line $x\geqslant0$, such that $F_{\eta}(x)>0$ for all $x>0$. 
 We claim that for all natural  $n<m/2-1$ the following equalities hold
\begin{align}\label{step3}
	F_{\eta_{m}+\eta}^{(n+1)}(0)=0\qquad \text{and}\qquad F_{\eta_{m}+\eta}^{(n)}(x)=o(F_{\eta_m+\eta}^{(n+1)}(x)),\quad x\to+0.
\end{align}
This means that~\eqref{step1} is still true after adding $\eta$.
According to step 1 it is sufficient to prove~\eqref{step0}  for $F_{\eta_{m}+\eta}(x)$ while $n<m/2-1$ . Note, that~\eqref{step0} holds for $F_{\eta_m}$.
Consider
\begin{align}
	\label{step3.1}
	F_{\eta_{m}+\eta}^{(n+2)}(x)
	=
	\frac{d^{n+2}}{dx^{n+2}}\int\limits_{-\infty}^{+\infty}F_{\eta_m}(x-y)\,dF_{\eta}(y)
	=
	\int\limits_{0}^{x}F_{\eta_m}^{(n+2)}(x-y)\,dF_{\eta}(y).
\end{align}
It is clear that $F_{\eta_{m}+\eta}^{(n+1)}(0)=0$. 
Due to mean value theorem there exists such $x_0\in(0,x)$ that
\begin{align*}
	F_{\eta_{m}+\eta}^{(n+2)}(x)
	=
	F_{\eta_m}^{(n+2)}(x_0)\int\limits_{0}^{x}\,dF_{\eta}(y)=F_{\eta_m}^{(n+2)}(x_0)(F_{\eta}(x)-F_{\eta}(0))>0.
\end{align*}

\textbf{\textit{Step 4:}} Let us prove that Lemma~\ref{lemmaF} is true for \emph{infinite} sum, that is for the function $F(x)$. 
Fix $n$. Choose $m\in\mathbb{N}$ such that $m>2(n+1)$. Write $\xi$ as the following sum
\begin{align*}
	\xi
	=
	\eta_m+\sum\limits_{p=m+1}^{\infty}\mu_p\xi_p^2
	=:
	\eta_m+\eta.
\end{align*}
It is clear that $\eta$ has a positive mass in any neighborhood of zero. Thus, using step 3 for $\xi$, we finish the proof.
\endproof

Let us prove the following lemma of Tauberian type (on differentiability of asymptotics).
\begin{lm}\label{Tauberian}
	Let $f(x):[0,+\infty)\to \mathbb{R}$,  $f'(x)$ is non-decreasing and 
	\begin{align}
	f(x)\sim \mathcal{A}\,x^{\alpha}L(x)\exp(-\mathcal{D}x^{-\beta}),\quad \text{ as }x\to+0,\qquad \alpha\in\mathbb{R}, C>0,\beta>0, \mathcal{D}>0,
	\end{align}
	where $L(x)>0$ is a slowly varying function at zero. 
	Then
	\begin{align}
	f'(x)\sim \mathcal{A}\,\mathcal{D}\beta x^{\alpha-\beta-1}L(x)\exp(-\mathcal{D}x^{-\beta}),\quad \text{ as }x\to+0.
	\end{align}
\end{lm}
\beginproofdot It is sufficient to prove the following inequalities
\begin{align}
\label{Tauberian_limsup}
\limsup\limits_{x\to+0}f'(x)/(x^{\alpha-\beta-1}L(x)\exp(-\mathcal{D}x^{-\beta}))\leqslant \mathcal{A}\,\mathcal{D}\beta;
\\
\label{Tauberian_liminf}
\liminf\limits_{x\to+0}f'(x)/(x^{\alpha-\beta-1}L(x)\exp(-\mathcal{D}x^{-\beta}))\geqslant \mathcal{A}\,\mathcal{D}\beta. 
\end{align}
Let us first prove \eqref{Tauberian_limsup}. By the mean value theorem there exists a point $\xi\in(x,x+h)$ such that $f(x+h)-f(x)=hf'(\xi)\geqslant h f'(x)$ for $h>0$. Fix $\varepsilon>0$. Then if $x$ is sufficiently small and $0<h<x$, we have
\begin{align*}
f'(x)&\leqslant \frac{\mathcal{A}\,(x+h)^{\alpha}L(x+h)\exp(-\mathcal{D}(x+h)^{-\beta})(1+\varepsilon)-\mathcal{A}\,x^{\alpha}L(x)\exp(-\mathcal{D}x^{-\beta})(1-\varepsilon)}{h}
\\
&=\frac{\mathcal{A}\,(x+h)^{\alpha}L(x+h)\exp(-\mathcal{D}(x+h)^{-\beta})-\mathcal{A}\,x^{\alpha}L(x+h)\exp(-\mathcal{D}(x+h)^{-\beta})}{h}
\\
&+\frac{\mathcal{A}\,x^{\alpha}L(x+h)\exp(-\mathcal{D}(x+h)^{-\beta})-\mathcal{A}\,x^{\alpha}L(x)\exp(-\mathcal{D}(x+h)^{-\beta})}{h}
\\
&+\frac{\mathcal{A}\,x^{\alpha}L(x)\exp(-\mathcal{D}(x+h)^{-\beta})-\mathcal{A}\,x^{\alpha}L(x)\exp(-\mathcal{D}x^{-\beta})}{h}+\frac{\varepsilon \mathcal{A}\,x^{\alpha}L(x)\exp(-\mathcal{D}x^{-\beta})}{h}\\
&+\frac{\varepsilon \mathcal{A}\,(x+h)^{\alpha}L(x+h)\exp(-\mathcal{D}(x+h)^{-\beta})}{h}=:B_1+B_2+B_3+B_4+B_5.
\end{align*}
Since $((x+h)^{\alpha}-x^{\alpha})/h\leqslant \alpha x^{\alpha-1}+Chx^{\alpha-2}$, we have
\begin{align*}
\frac{B_1}{x^{\alpha-\beta-1}L(x)\exp(-\mathcal{D}(x+h)^{-\beta})}&\leqslant (\mathcal{A}\,\alpha x^{\beta}+h\mathcal{A}\,Cx^{\beta-1})\frac{L(x+h)}{L(x)}.
\end{align*}
Without loss of generality we may assume that 
\begin{align*}
L(x)=\exp\left(\int\limits_B^x\frac{\eta(y)}{y}\,dy\right),
\end{align*}
where $\eta(y)$ is a measurable function on $[B,\infty)$ and $\eta(y)\to0$ as $y\to\infty$. 
This means that 
\begin{align}\label{SVF}
\frac{xL'(x)}{L(x)}\to0.
\end{align}
Further,
\begin{align}\nonumber
\frac{B_2}{x^{\alpha-\beta-1}L(x)\exp(-\mathcal{D}(x+h)^{-\beta})}&= \mathcal{A}\cdot\frac{x^{\beta+1}}{L(x)}\cdot\frac{\left(L(x+h)-L(x)\right)}{h}=\mathcal{A}\cdot\frac{x^{\beta+1}L'(x+z)}{L(x)}
\end{align}
for some $z\in(0,h)$. 
Next,
\begin{align*}
\frac{B_3}{x^{\alpha}L(x)\exp(-\mathcal{D}x^{-\beta})}&= \frac{\mathcal{A}}h\left(\exp(-\mathcal{D}(x+h)^{-\beta}+\mathcal{D}x^{-\beta})-1\right)=\frac{\mathcal{A}}h\left(\exp\left(\mathcal{D}\frac{(1+h/x)^{\beta}-1}{(x+h)^{\beta}}\right)-1\right)
\\
&\leqslant \frac{\mathcal{A}}h\left(\exp\left(\mathcal{D}\frac{(1+h/x)^{\beta}-1}{x^{\beta}}\right)-1\right).
\end{align*}
Using inequalities $(1+y)^\beta-1\leqslant y\beta+Cy^2$ and $\exp(y)-1\leqslant y+Cy^2$ for small enough $y$, we get
\begin{align*}
\frac{B_3}{x^{\alpha}L(x)\exp(-\mathcal{D}x^{-\beta})}&\leqslant \frac{\mathcal{A}}h\left(\exp\left(\mathcal{D}\frac{\beta h/x+C(h/x)^2}{x^{\beta}}\right)-1\right)
\\
&\leqslant
\frac{\mathcal{A}}h\left(\mathcal{D}\frac{\beta h/x+C(h/x)^2}{x^{\beta}}+C\frac{(h/x)^2}{x^{2\beta}}\right)
\\
&=
\mathcal{A}\,\mathcal{D}\beta x^{-\beta-1}+\mathcal{A}\,Ch\mathcal{D} x^{-\beta-2}+\mathcal{A}\,Ch x^{-2\beta-2}.
\end{align*}
In a similar way we obtain
\begin{align*}
\frac{B_5}{x^{\alpha}L(x)\exp(-\mathcal{D}x^{-\beta})}&\leqslant \frac{\varepsilon \mathcal{A}}{h}\frac{L(x+h)}{L(x)}\exp\left(\mathcal{D}\frac{\beta h/x+C(h/x)^2}{x^{\beta}}\right).
\end{align*}
Let $h=\sqrt{\varepsilon}x^{\beta+1}$. Then
\begin{align*}
B_1&\leqslant x^{\alpha-\beta-1}L(x)\exp(-\mathcal{D}(x+h)^{-\beta})\cdot ( \mathcal{A}\,\alpha x^{\beta}+\sqrt{\varepsilon} \mathcal{A}\,Cx^{2\beta})\cdot \frac{L(x+\sqrt{\varepsilon}x^{\beta+1})}{L(x)};
\\
B_2&\leqslant x^{\alpha-\beta-1}L(x)\exp(-\mathcal{D}(x+h)^{-\beta})\cdot \mathcal{A}\cdot\frac{x^{\beta+1}L'(x+z)}{L(x)},\qquad z\in(0,\sqrt{\varepsilon}x^{\beta+1});
\\
B_3&\leqslant x^{\alpha-\beta-1}L(x)\exp(-\mathcal{D}x^{-\beta})\cdot\left(\mathcal{A}\,\mathcal{D}\beta+\sqrt{\varepsilon}\mathcal{A}\,C\mathcal{D} x^{\beta}+\sqrt{\varepsilon}\mathcal{A}\,C\right);
\\
B_4&\leqslant x^{\alpha-\beta-1}L(x)\exp(-\mathcal{D}x^{-\beta})\cdot \sqrt{\varepsilon}\mathcal{A};
\\
B_5&\leqslant x^{\alpha-\beta-1}L(x)\exp(-\mathcal{D}x^{-\beta})\cdot\sqrt{\varepsilon}\mathcal{A}\,\exp(\mathcal{D}\beta\sqrt{\varepsilon}+\mathcal{D}Cx^{\beta}\varepsilon)\cdot\frac{L(x+\sqrt{\varepsilon}x^{\beta+1})}{L(x)}.
\end{align*}
Let us note that as $L(x)$ is a slowly varying function we obtain
\begin{align*}
\lim\limits_{x\to0}\frac{L(x+\sqrt{\varepsilon}x^{\beta+1})}{L(x)}=1.
\end{align*}
Indeed,
\begin{align*}
\lim\limits_{x\to0}\frac{L(x+\sqrt{\varepsilon}x^{\beta+1})}{L(x)}&=1+\lim\limits_{x\to0}\frac{L(x+\sqrt{\varepsilon}x^{\beta+1})-L(x)}{L(x)}=1+\lim\limits_{x\to0}\frac{\sqrt{\varepsilon}x^{\beta+1}L'(x+z)}{L(x)}
\\
&=1+\sqrt{\varepsilon}\lim\limits_{x\to0}\frac{x^{\beta/2}L(x+z)}{x^{-\beta/2}L(x)}\cdot\frac{x}{x+z}\cdot\frac{(x+z)L'(x+z)}{L(x+z)}=1.
\end{align*}
for some $z\in(0,\sqrt{\varepsilon}x^{\beta+1})$. The last equality holds due to \eqref{SVF} and the fact that $x^{\alpha}L(x)\to 0$ as $x\to0$ for any $\alpha>0$.

From these estimates it follows that $$
\frac{f'(x)}{x^{\alpha-\beta-1}L(x)\exp(-\mathcal{D}x^{-\beta})}\leqslant \mathcal{A}\,\mathcal{D}\beta+C(\sqrt{\varepsilon}+x^{\beta}).$$
This completes the proof of \eqref{Tauberian_limsup}. To obtain \eqref{Tauberian_liminf} we proceed in a similar way using the inequality $f(x)-f(x-h)\leqslant hf'(x)$.
\endproof

\beginproof of Lemma \ref{diff}.
Combining Lemmata \ref{lemmaF} and  \ref{Tauberian} we immediately get Lemma \ref{diff}.
\endproof

Petrova Yulia: yu.pe.petrova@yandex.ru

Chebyshev Laboratory

St. Petersburg State University, 14th Line V.O., 29B, Saint Petersburg 199178 Russia.


\end{document}